\documentclass[12pt]{amsart}
\usepackage{amsmath}
\usepackage{amssymb}
\usepackage{amscd}
\usepackage[all, cmtip]{xy}
\def\lra{\longrightarrow}
\def\map#1{\,{\buildrel #1 \over \lra}\,}
\def\smap#1{{\buildrel #1 \over \to}}

\newcommand{\pc}{^{\mathrm{pc}}}

\newcommand{\into}{\rightarrowtail}
\newcommand{\onto}{\twoheadrightarrow}

\def\div{{\mathrm{div}}}
\def\bM{\mathbf M}
\def\cA{\mathcal A}

\def\cC{\mathcal C}

\def\cM{\mathcal M}
\def\cN{\mathcal N}
\def\cS{\mathcal S}

\def\T{\mathbb T}
\def\inv{^{-1}}
\def\ie{i.e.,\ }
\DeclareMathOperator{\id}{id}

\newcommand{\Z}{\mathbb{Z}}
\newcommand{\N}{\mathbb{N}}
\newcommand{\F}{\mathbb{F}}
\newcommand{\bS}{\mathbb{S}}
\newcommand{\X}{S} 
\newcommand{\m}{\mathfrak{m}}
\newcommand{\G}{{\Gamma}}
\newcommand{\colim}{\operatorname*{colim}}
\newcommand{\coker}{\operatorname{coker}}
\newcommand{\Hom}{\operatorname{Hom}}
\newcommand{\Aut}{\operatorname{Aut}}
\newcommand{\Cl}{\operatorname{Cl}}
\newcommand{\Pic}{\operatorname{Pic}}
\def\MSpec{\operatorname{MSpec}}

\newcommand{\Sets}{\mathbf{-Sets}}
\newcommand{\Div}{\mathrm{Div}}
\def\fin{{\mathrm{fin}}}

\def\im{{\mathrm{im}}}

\DeclareMathOperator{\height}{ht}
\def\double{^{[d]}}

\newcommand{\mono}{\ar@{>->}}
\newcommand{\epi}{\ar@{->>}}
\newcommand{\emap}{\ar@{o->}}

\numberwithin{equation}{section}

\theoremstyle{plain}
\newtheorem{thm}[equation]{Theorem}
\newtheorem{prop}[equation]{Proposition}
\newtheorem{lem}[equation]{Lemma}

\newtheorem{substuff}{Remark}[equation]

\theoremstyle{definition}
\newtheorem{defn}[equation]{Definition}
\newtheorem{ex}[equation]{Example}
\newtheorem{subex}[substuff]{Example}

\newtheorem{app}[equation]{Application}
\newtheorem{construct}[equation]{Construction}

\theoremstyle{remark}
\newtheorem{rem}[equation]{Remark}
\newtheorem{subrem}[substuff]{Remark}

\begin {document}
\title{Localization, monoid sets and $K$-theory}
\date{\today}
\author{Ian Coley}
\address{Math.\ Dept., Rutgers University, New Brunswick, NJ 08901, USA}
\email{iacoley@math.rutgers.edu}
\urladdr{http://iancoley.org}

\author{Charles Weibel}
\address{Math.\ Dept., Rutgers University, New Brunswick, NJ 08901, USA}
\email{weibel@math.rutgers.edu}\urladdr{http://math.rutgers.edu/$\sim$weibel}
\thanks{Weibel was supported by NSF grants, and the Simonyi Endowment at IAS}

\begin{abstract}
We develop the $K$-theory of sets with an action of a pointed monoid
(or monoid scheme), analogous to the $K$-theory of modules over a 
ring (or scheme).
In order to form localization sequences, we construct the quotient
category of a nice regular category by a Serre subcategory.
\end{abstract}
\maketitle
There are several interesting categories of $A$--sets when 
$A$ is a noetherian pointed monoid $A$.
One is the category of all finitely generated $A$--sets; 
it is studied in \cite{Deitmar}, \cite{CLS} and \cite{FW}.
Another is the category of all 
{\it partially cancellative} $A$--sets, which
 is studied in \cite{HW}; see Definition \ref{def:pc} below.
These $A$--sets are well-behaved, include free $A$--sets, and 
turn out to be useful in toric geometry; see \cite{CHWW}.
There are similar categories over a noetherian monoid scheme,
including categories of $A$--sets with a restriction 
on the codimension of support.

This paper is primarily concerned with the $K$-theory
of these categories. One of our themes is that they closely track
Quillen's $K'$-theory of noetherian rings and schemes.
In particular, we show that their $K$-theory has localization sequences
and coniveau spectral sequences. In particular,
we prove in Theorem \ref{thm:Gersten} that,
if $\X$ is a 0-smooth monoid scheme, and $\mathcal{K'}_*$ 
is the sheaf associated to the $K$-theory of
partially cancellative $\cA_\X$--sets
then 
\[
E_2^{p,q} = H^p(\X,\mathcal{K}'_{-q}).
\]

The fundamental new tool we develop and use is the construction of
the localization $\cM/\cC$ of a nice regular category $\cM$
(such as the above categories of $A$--sets) by a
``Serre subcategory'' $\cC$. This construction generalizes 
Gabriel's localization of abelian categories, and is of independent interest.
For expositional reasons, the construction of $\cM/\cC$ is
postponed until Sections \ref{sec:cat-thy} and \ref{sec:compare};
we summarize it in Section \ref{sec:local}, 
where we apply it to obtain $K$-theory localization sequences
\[ K(\cC) \to K(\cM) \to K(\cM/\cC). \]

\smallskip
This paper organized as follows. In the first few sections, we
recall the basic definitions of $A$--sets, 
quasi-exact categories (such as $A$--sets) and their $K$-theory.
Section \ref{sec:cgw} recalls the Campbell--Zakharevich notions 
of {\it CGW} and {\it ACGW} categories
which are needed for $K$-theory localization sequences.
Our version of the localization sequence is summarized in
Section \ref{sec:local}, with the technical details relegated to
Sections \ref{sec:cat-thy} and \ref{sec:compare}.  
Sections \ref{sec:FT}--\ref{sec:GC}
construct the analogues of the Gersten--Quillen spectral
sequence associated to the coniveau filtration, and 
establish an analogue of Quillen's proof of Gersten's Conjecture.

\subsection*{Acknowledgements}
The authors would like to thank Inna Zakharevich for her help in
understanding CGW-categories in \cite{CZ}.

\section{monoids and monoid sets}\label{sec:pc}

By an (abelian) {\it monoid} $A$ we mean a pointed set with a commutative
product and distinguished elements $\{0,1\}$ such that
$a\cdot1=a$ and $a\cdot0=0$ for all $a\in A$; these
are sometimes called ``pointed monoids.''
For example, the free monoid $\N$ is 
$\{0,1,t,t^2,...\}$, and if $A$ and $B$ are monoids, so is $A\wedge B$.
The initial pointed monoid $\{0,1\}$ is called $\mathbb{F}_1$ in \cite{CLS}.
We call the set $A^\times$ of invertible elements of a monoid its
\emph{units}; these form a group.

The monoid spectrum $\MSpec(A)$ is the space of prime ideals in $A$,
together with the sheaf whose stalk at a prime $s$ 
is the localization $A_s$ of the monoid $A$.
A {\em monoid scheme} $(\X,\cA_\X)$ is a topological space $\X$ with
a sheaf of monoids $\cA_\X$ which is locally $\MSpec$ of a monoid.

If $A$ is a (pointed) monoid,
an {\em $A$--set} is a pointed set $(X,\ast)$ with an action of $A$;
in particular, $1\cdot x=x$ and $0\cdot x=\ast$ for all $x\in X$.
For example, a free $A$--set is just a wedge of copies of $A$.
We write $A\Sets$ for the evident category of $A$--sets.

If $(\X,\cA_\X)$ is a monoid scheme, we write 
$\cA_\X\Sets$ for the category of (quasi-coherent) sheaves of $\cA_\X$--sets.
If $\X=\MSpec(A)$, $\cA_\X\Sets$ is equivalent to $A\Sets$.

A monoid $A$ is {\em partially cancellative}, or {\em pc} for short, if 
$ac=bc\ne\ast$ implies $a=b$ for all $a,b,c$ in $A$. 
If $A$ is a cancellative monoid then $A$ and $A/I$ are pc
for every ideal $I$ of $A$; 
the prototypical finite pc monoid is $\{1,t,...,t^N=t^{N+1}=\ast\}.$

\begin{defn}\label{def:pc}
If $A$ is a pc monoid, we say that a pointed $A$--set $X$ is 
{\em partially cancellative} if for every $x\in X$ and $a,b$ in $A$, 
if $a\cdot x=b\cdot x\ne\ast$ then $a=b$.
The subcategory of pc monoid sets is 
closed in $\cA\Sets$ under subobjects and quotients.
We do not impose the condition that $a\cdot x = a\cdot y\ne\ast$ implies
$x=y$, as this rules out some key examples (see Example \ref{ex:rooted}).

Similarly, we say a monoid scheme $(\X,\cA_\X)$ is \emph{pc} if it is locally
pc, and a sheaf of $\cA_\X$--sets is pc if it is locally pc.
If $A$ is pc, open subschemes of $\MSpec(A)$ are pc
but are not usually affine.
\end{defn}

\begin{ex}\label{ex:rooted}
A (pointed) $\N$-set is just a pointed set $X$ with a successor
function $x\mapsto tx$. 
Every finite rooted tree is a pc $\N$-set;
the successor of $x$ is the adjacent vertex closer to the root vertex $\ast$.
In fact, a finite $\N$-set is partially cancellative if and only
if it is a rooted tree, because for every $x\in X$, the sequence 
$\{x,tx,t^2x,...\}$ terminates at the root vertex.

An $\N$-set is pc if and only if it contains no loop, \ie
for each element $x\ne\ast$ and integer $d>1$ we have $t^dx\ne x$.
A typical non-pc set is 
$\{1,t,...,t^d,...,t^N,\ast:t^{N}=t^{d}\}$.
\end{ex}

\begin{ex}\label{G-sets}
Let $A$ be $\G_+$, where $\G$ is an abelian group, and 
`$+$' is a disjoint basepoint. Then every $A$--set is a 
wedge of copies of cosets $(\G/H)_+$.
If $H$ is a proper subgroup of $\G$, then $(\G/H)_+$ is 
not pc because $h\cdot x=1\cdot x$ for $h\in H$. 
Thus a $\G_+$-set is pc if and only if it is free.
\end{ex}

\section{$K$-theory of quasi-exact categories}\label{sec:quasi}

The notion of exact category \cite{Q} has the following generalization.

\begin{defn}(\cite{Deitmar}, \cite[Ex.\,IV.6.14]{WK})
A {\it quasi-exact category} is a category $\cM$
with a distinguished zero object,
and a coproduct $\vee$, equipped with a family $\cS$ of sequences of the form
\begin{equation}\label{eq:admissible}
X \map{i} Y  \map{p} Z, 
\end{equation}
called ``admissible,'' such that: \\
(i) any sequence isomorphic to an admissible sequence is admissible;
\\
(ii) for any admissible sequence \eqref{eq:admissible},
$p$ is a cokernel for $i$ and $i$ is a kernel for $p$;\\
(iii) $\cS$ contains all split sequences (those with $X\cong Y\vee Z$);
and \\(iv) the class of admissible epimorphisms 
(resp., admissible monics) is closed under composition and
pullback along admissible monics (resp., pullback along admissible
epimorphisms). 

We will write admissible monics as $X\into Y$ and 
admissible epics as $Y\onto Z$, and will
often write $Y/X$ for the cokernel of $X\into Y$.
\end{defn}

The prototype of a quasi-exact category is a {\it regular}
category; see Definition \ref{def:regular}. The exact
sequences are the sequences \eqref{eq:admissible}
for which $p$ is a cokernel for $i$ and $i$ is a kernel for $p$.

\begin{subex}
The category $\mathbf{Sets}$ of sets is quasi-exact.
More generally,
if $A$ is a (pointed) monoid, the category $A\Sets$ is quasi-exact;
a sequence \eqref{eq:admissible} is admissible if $X\into Y$ is an
injection, and $Z$ is isomorphic to the quotient $A$--set $Y/X$.
If $A$ is a noetherian monoid, the category $\cM(A)$
of finitely generated pointed $A$--sets is quasi-exact.
(See \cite[Ex.\,IV.6.16]{WK}.)

Similarly, if $(\X,\cA_\X)$ is a monoid scheme,
$\cA_\X\Sets$ is  a quasi-exact category;
the admissible sequences are defined locally.
\end{subex}

The wedge $X\wedge_{\cA} Y$ of two $\cA$--sets is defined to be
the quotient of $X\wedge Y$ by the equivalence relation that
$(xa,y)\sim(x,ay)$ for all $a\in A$ $x\in X$ and $y\in Y$.
As noted in \cite[5.10]{CLS}, $X\wedge_{\cA} Y$ is exact
in both variables; 
in fact, $\cA$--Sets is a symmetric monoidal category.

\begin{defn}
If $\cM$ is quasi-exact,
Quillen's construction in \cite{Q} yields a category $Q\cM$, 
and $K(\cM)$ is the connective spectrum with initial space
$\Omega BQ\cM$; we write $K_n(\cM)$ for $\pi_nK(\cM)$.
The group $K_0(\cM)$ is generated by the objects, modulo the relations
that $[X]=[Y]+[Z]$ for every sequence \eqref{eq:admissible}.
\end{defn}

\begin{subex}
The category $\mathbf{Sets}_\fin$ of finite pointed sets
is quasi-exact; every admissible sequence is split exact.
It is well known that the Barratt--Priddy--Quillen theorem implies that
$K(\mathbf{Sets}_\fin)$ is homotopy equivalent to
the sphere spectrum $\bS$.
(See \cite{Deitmar}, \cite{CLS}, \cite[Ex.\,IV.6.15]{WK}.)
\end{subex}

\begin{defn}\label{def:G(A)}
If $A$ is noetherian, the category $\cM(A)$ of noetherian $A$--sets
is quasi-exact; following \cite{CLS}, we write $G(A)$ for $K(\cM(A))$.
Similarly, if $\X$ is a noetherian monoid scheme,
the category $\cM(\X)$ of sheaves of noetherian $\cA_\X$--sets
is quasi-exact, and we write $G(\X)$ for $K(\cM(\X))$. 
Since $\wedge_\cA$ is biexact, $G(\X)$ is an $E_\infty$-ring spectrum.
See \cite[IV.6.6]{WK} and \cite[3.8.2]{Bar}.

$G$-theory is contravariantly functorial with 
respect to monoid maps $A\to B$ for which $B$ is 
noetherian as an $A$-module, and covariant for
flat maps (a monoid map
$A\to B$ is \emph{flat} if the base extension 
$X\mapsto B\wedge_A X$ is exact).

If $A$ is a noetherian pc monoid, the category $\cM\pc(A)$ of 
noetherian pc $A$--sets is quasi-exact;
following \cite{HW}, we write $K'(A)$ for $K(\cM\pc(A))$.
Similarly, if $(\X,\cA_\X)$ is a noetherian pc monoid scheme, 
the category $\cM\pc(\X)$ of sheaves of noetherian pc $\cA_\X$--sets is
quasi-exact, and we write $K'(\X)$ for $K(\cM\pc(\X))$. 
$\cM\pc(\X)$ is symmetric monoidal,
since $X\wedge_\cA Y$ is pc when both $X$ and $Y$ are, and
$K'(A)$ and $K'(\X)$ are $E_\infty$-ring spectra.

As noted in \cite[3.1.1]{HW},
$K'$ is covariantly functorial for flat maps.
\end{defn}

\begin{subex}\label{K(G)}
Let $\G$ be a group; the $K$-theory of the category of 
finitely generated free $\G$--sets is written as $K(\G)$.  
When $A=\G_+$, we saw in Example \ref{G-sets} that
every pc $A$--set is a free $A$--set, so $K'(\G_+)\simeq K(\G)$.
If $X=A^{\vee r}$, then $\Aut(X)$ is the wreath product $\G\wr\Sigma_r$.
By the Barratt--Priddy--Quillen theorem,
\[ 
K'(\G_+) \simeq K(\G) \simeq S^\infty(B\G_+) = S^\infty(B\G)\vee\bS.
\]

Similarly, if $\G$ is abelian, the category of $\G$--sets is the 
product over all subgroups $[H]$ of the category of $\G/H$--sets%
\footnote{If $\G$ is not abelian, replace $\G/H$ by the
Weyl group $N_\Gamma(H)/H$}. 
By the equivariant Barratt--Priddy--Quillen theorem \cite[5.1]{CDD},
\[
G(\G_+) \simeq \bigvee\nolimits_{H\le\G} S^\infty[B(\G/H)_+].
\]
In particular, $G_0(\G_+)=\pi_0 G(\G_+)$ is the Burnside ring of $\G$.
\end{subex}

We say that an $A$--set $X$ has 
\emph{finite length} if it has a finite filtration
\[ 1\into F_1X\into\cdots\into F_nX=X \]
such that each $F_iX/F_{i-1}X$ is irreducible,
\ie isomorphic to $(A^\times)_+$. If $A$ is noetherian, 
the category $\cM_\fin(A)$
of $A$--sets of finite length is quasi-exact,
as is the category $\cM\pc_\fin(A)$ of pc
$A$--sets of finite length. 
By D\'evissage (see \cite[6.2]{CZ} or \cite[3.3, 3.3.1]{HW}), we have:

\begin{thm}\label{thm:devissage}
If $A$ is a noetherian monoid, with units $\G=A^\times$, then 
$K(\cM_\fin(A))\simeq G(\G_+)$ and if $A$ is a pc monoid then:
\[ 
K(\cM\pc_\fin(A))\cong K(\G)
\simeq S^\infty(B\G_+)\simeq S^\infty(B\G)\vee\bS.
\]
If $A$ has finite length 
then  $G(A)\simeq G(\G_+)$ and $K'(A)\simeq K(\G).$
\end{thm}

The following result is taken from \cite[5.3]{HW}.

\begin{thm}\label{thm:X-U}
  Let $\X$ be a pc monoid scheme and $Z\map{i} \X$ 
  an equivariant closed subscheme with open complement $U\map{j} \X$. 
  Then there is a fibration sequence of spectra
  \[
K'(Z)\map{i_*} K'(\X) \map{j^*} K'(U).
\]
\end{thm}

\begin{rem}
By \cite[1.3]{FW}, every ideal $I$ in a noetherian monoid 
has a finite number of associated primes. It follows that
the support of an equivariant closed subscheme $Z$ of a noetherian monoid
scheme has a finite number of minimal points.
\end{rem}

\goodbreak
\section{CGW and ACGW categories}\label{sec:cgw}

Campbell and Zakharevich have defined a \emph{CGW-category} to be a double
category satisfying a certain list of axioms; we refer
the reader to \cite[2.5]{CZ} for the precise definition.
Here is our main example:

If $\cM$ is a quasi-exact category, we can form a double category
$\cM\double$ with the same objects as $\cM$; the horizontal and vertical
maps are the admissible monics and epics (composed backwards), respectively,
and the 2--cells are commutative diagrams of the form
\begin{equation}\label{ME-square}
\vcenter{\xymatrix@R=1.9em{
X~ \ar@{>->}[r] & Y \\
X'~ \ar@{>->}[r]\ar@{->>}[u] & Y'\ar@{->>}[u]
}}
\end{equation}
We say that a square \eqref{ME-square}
is {\it distinguished} if the natural map of cokernels
$Y'/X'\to Y/X$ is an isomorphism.  Thus distinguished squares are both 
pushout squares and pullback squares.

As pointed out in \cite[2.3]{HW}, 
if $\cM$ is a quasi-exact category, $\cM\double$
is a CGW-category, and $\cM$ is an ``ambient category'' 
for $\cM\double$ in the sense of \cite[2.3]{CZ}, where
$k(Y\onto Y/X)$ is its kernel $X\into Y$ and 
$c(X\into Y)$ is its cokernel $Y\onto Y/X$.
Moreover, the $K$-theory $K(\cM\double)$ agrees with $K(\cM)$.

There is a stronger notion, that of an ACGW-category, defined in
\cite[5.6]{CZ}. If $\cM$ is an ACGW-category and $\cC$ is an ACGW-subcategory
closed under subobjects, quotients, and extensions 
(see \cite[2.12]{CZ}), the double category $\cM\backslash \cC$
is defined in \cite[8.1]{CZ}. There is a canonical morphism of double categories
$\cM\to\cM\backslash \cC$.

A key construction in the language of ACGW-categories is 
\cite[8.6]{CZ}, which we now cite:

\begin{thm}\label{CZ8.6}
Let $\cM$ be an ACGW-category and let $\cC$ be an ACGW-subcategory 
satisfying the technical conditions (W), (E) and (CGW). 
Then 
\begin{equation*}
K(\cC)\to K(\cM)\to K(\cM\backslash \cC)
\end{equation*}
is a homotopy fiber sequence.
\end{thm}

The condition (CGW) is that $\cM\backslash \cC$ is a CGW category; 
condition (W) is described in Section~\ref{sec:compare}.
Condition (E) is unnecessary; see Remark \ref{rem:(E)}

\begin{ex}
By \cite[2.5]{HW}, if $\cM$ is a quasi-exact subcategory of
$A\Sets$ (or $\cA_\X\Sets$) closed under 
pushouts along pairs of monics, and pullbacks along pairs of epics,
then the associated double category 
$\cM\double$ is an ACGW-category.

In particular, when $A$ is a noetherian monoid, 
the associated double categories of $\cM(A)$ 
and $\cM\pc(A)$ are ACGW-categories;
the same is true for the subcategories
$\cM(\X)$ and $\cM\pc(\X)$ over a monoid scheme.
\end{ex}

The difficulty with Theorem \ref{CZ8.6} is that $\cM\backslash \cC$
is hard to work with.  Since we are working with 
ambient categories having more structure,
we have access to a different kind of localization, viz.,
a quotient category $\cM/\cC$,
which we establish in Theorem \ref{thm:M/C-long}, 
and summarize in the next section.

\goodbreak
\section{The localization $\cM/\cC$}\label{sec:local}

To define $\cM/\cC$, we need to introduce
some category-theoretic vocabulary. We say that
a pointed category is \emph{regular} if it admits all finite limits 
and has a good intrinsic notion of short exact sequence, 
where the admissible monics are always kernels and the admissible
epics are always cokernels. (See Definition \ref{def:regular}.)
A \emph{regular functor} is one that
preserves this structure, \ie finite limits and short exact sequences.

Every regular category is quasi-exact; see Remark \ref{reg=qexact}.
The categories $A\Sets$ and $\cM(A)$ are regular, 
as are  $A\Sets\pc$ and $\cM\pc(A)$. 

\begin{defn}
A full pointed subcategory $\cC$ of a regular category $\cM$
is a \emph{Serre subcategory} if it is closed under
finite limits and for every exact sequence 
\eqref{eq:admissible}, $Y$ is in $\cC$ if and only if 
both $X$ and $Z$ are in $\cC$.
\end{defn}

\begin{prop}\label{prop:Serre}
If $F\colon\cM\to\cN$ is a regular functor between regular categories,
then the full subcategory $\ker(F)=F^{-1}(0)$ is a Serre subcategory
of $\cM$.
\end{prop}

\begin{proof}
Clearly, $\ker(F)$ is closed under finite limits. Consider 
an exact sequence \eqref{eq:admissible} in $\cM$; applying $F$ yields
an exact sequence in $\cN$,
\begin{equation*}
F(X)\into F(Y)\onto F(Z).
\end{equation*}
It is clear
that $F(Y)\cong0$ if and only if $F(X)\cong F(Z)\cong0$.
\end{proof}

\begin{ex}\label{ex:S-torsion}
Let $S$ be a multiplicatively closed subset of an (abelian) monoid $A$.
Then the localization $A\Sets\to S^{-1}A\Sets$ is a regular functor, 
because it preserves exact sequences and finite limits.
Indeed, the ring-theoretic proofs go through; for example,
if $X\map{i}Y$ is an injection and $i(x)/s=0$ in $S^{-1}Y$
for $x\in X$ then some $i(s_1x)=0$ in $Y$ and hence $s_1x=0$ in $X$.

By Proposition \ref{prop:Serre}, the category of $S$-torsion $A$--sets is a
Serre subcategory of $A\Sets$; this description restricts to
Serre subcategories of $\cM(A)$ and $\cM\pc(A)$ 
(see Definition \ref{def:G(A)}).

More generally, if $U$ is open in a monoid scheme ($\X,\cA_\X)$ then
the localization $\cA_\X\Sets\to \cA_U\Sets$ is a regular functor.
Indeed, exact sequences and finite limits are determined locally on
affine open subschemes. 
By Proposition \ref{prop:Serre}, the category of $\cA_\X$--sets 
vanishing on $U$ is a Serre subcategory of $\cA_\X\Sets$, and this restricts
to Serre subcategories of $\cM(\X)$ and $\cM\pc(\X)$ as well. 
\end{ex}

Here are the two theorems about $\cM/\cC$ which we need.
The definition of an {\it adherent} category is given in 
Definition \ref{def:adherent};
it is a regular category satisfying some technical axioms.

\begin{thm}\label{thm:M/C}
Let $\cC$ be a Serre subcategory of an adherent category $\cM$. 
Then there exists a quotient category $\cM/\cC$, which is regular,
and a regular functor $\cM\to\cM/\cC$, which 
is initial among all regular functors from $\cM$ which send $\cC$ to zero.
\end{thm}

The construction is essentially a non-additive
version of Gabriel's construction \cite[Chapitre~III]{G}.
The full statement of Theorem \ref{thm:M/C}, and the proof, 
can be found in Theorem \ref{thm:M/C-long}. 

In Lemma \ref{lem:adherentACGW}, 
we point out that the double category $\cM^{[d]}$
associated with an adherent category $\cM$
is an ACGW category. 
A fuller statement of the following theorem (and its proof) is given in 
Lemma \ref{lem:M/C[d]} and Theorem \ref{thm:MC=MC-long}.

\begin{thm}\label{thm:MC=MC}
Let $\cC$ be a Serre subcategory of an adherent category $\cM$.
Then the double category associated to $\cM/\cC$ is a CGW-category,
equivalent to the CGW-localization $\cM\backslash\cC$ of 
$\cM\double$ by $\cC\double$, and there is a fibration sequence:
\[ K(\cC) \to K(\cM) \to K(\cM/\cC). \]
\end{thm}

\begin{ex}
Let $(\X,\cA_\X)$ be a monoid scheme.
For brevity, we write $\bM(\X)$ for the category $\cA_\X\Sets$.
If $\X=\MSpec(A)$, we will write $\bM(A)$ for $\bM(\X)$.
If  $U\subset \X$ is an open subscheme, with complement $Z$,
then $\bM(U)$ is $\bM(\X)/\bM_Z(\X)$, where
$\bM_Z(\X)$ is the subcategory of all $M$ in
$\bM(\X)$ with $M|_U=0$. Then
$\bM(U)$ is $\bM(\X)/\bM_Z(\X)$.
The proof is the same as Gabriel's in \cite[p.\,380]{G};
use the fact that the exact functor $j^*\colon\bM(\X) \to\bM(U)$ has a right
adjoint, 
the \emph{section functor} $j_*\colon\bM(U)\to\bM(\X)$, and the universal
map $j_*\colon\bM(U)\to \bM(\X)/\bM_Z(\X)$ is an equivalence.

Suppose $\X$ is noetherian.
As in Gabriel \cite[p.\,379]{G}, there is a bijection between
the Serre subcategories of $\bM(\X)$ and $\cM(\X)$ 
(noetherian $\cA_\X$-sets), given by 
$\bM'(\X) \mapsto \cM'(\X)=\bM'(\X)\cap\cM(\X)$.
Hence $\cM_Z(\X)$ is a Serre subcategory of $\cM(\X)$.
By Theorem \ref{thm:MC=MC} and D\'evissage \cite[6.2]{CZ},
there is a fibration sequence
\[
G(Z) \to G(\X) \to G(U).
\]
Similarly, if $\X$ is pc, then
$\cM\pc_Z(\X)$ is a Serre subcategory of $\cM\pc(\X)$,
and we recover the fibration sequence of Theorem \ref{thm:X-U}.
\end{ex}

\begin{app}\label{app:M/M1}
If $(\X,\cA_\X)$ is a monoid scheme, then
the subcategory $\bM_d(\X)$ of sheaves in $\bM(\X)$ which are
supported in codimension $\ge d$ is a Serre subcategory.
It contains $\bM_{d+1}(\X)$ and
\begin{equation*}
\bM_d(\X) / \bM_{d+1}(\X)\cong\bigvee\nolimits_{\height(s)=d} \bM_{s}(\cA_s),
\end{equation*}
where $\bM_{s}(\cA_s)$ is the category of torsion $\cA_s$--sets.
Similarly, if $\X$ is noetherian then
$\cM_d(\X) / \cM_{d+1}(\X)$ is equivalent to the wedge of
the categories $\cM_\fin(\cA_s)$ of
finite length $\cA_s$--sets. The parallel assertion
holds for $\cM\pc_d(\X) / \cM\pc_{d+1}(\X)$.

For example, if $\X$ is noetherian, then $\bM_1(\X)$
is the category of $\cA_\X$--sets vanishing on the open subscheme
$U$ consisting of the (finitely many) points of codimension 0.
In this case, we get the equivalences
\begin{equation*}
\bM(\X)/\bM_1(\X) \cong \bM(U)=\bigvee\nolimits_{\height(s)=0} \bM(\cA_s)
\end{equation*}
and $\cM(\X)/\cM_1(\X) \cong \cM(U)=\bigvee_{\height(s)=0} \cM(\cA_s).$

If $\X$ has a unique minimal point,  $j_*$
identifies $\bM(\X)/\bM_1(\X)$ with the category $Q$ of 
torsionfree objects in $\bM(\X)$,
and identifies $\cM(\X)/\cM_1(\X)$ with the noetherian objects in $Q$.
\end{app}

\section{The coniveau filtration}\label{sec:FT}

We will apply the localization theorem \ref{thm:MC=MC}
to $\cM_{d+1}(\X) \subset \cM_d(\X)$ and its pc analogue.
By \ref{app:M/M1},
$\cM_d\pc(\X) / \cM_{d+1}\pc(\X)$
 is equivalent to the coproduct over height $d$ points $s$
of torsion pc $\cA_s$-modules.

Let $(\X,\cA_\X)$ be a noetherian monoid scheme, such as $\MSpec(A)$.
If $s$ is a point in $\X$, we write
$\G_s$ for the group of units of the pointed monoid $\cA_s$, so
$(\G_s)_+=\cA_s/p_s$. If $\X=\MSpec(A)$, then
$p_s$ is a prime ideal in the monoid $A$, and $\G_s$ is 
the group completion of the cancellative monoid $A/p_s$. 
By Example \ref{K(G)}, $K'(\G_s)= S^\infty(B\G_s)\vee \bS$.

Let $\X$ be a noetherian monoid scheme. 
By Theorem \ref{thm:MC=MC}, the quasi-exact category 
$\cM_d(\X)/\cM_{d+1}(\X)$ is the ambient category of
the double category $\cM_d(\X) \backslash \cM_{d+1}(\X)$.
By \ref{app:M/M1}, the quotient $\cM_d(\X)/\cM_{d+1}(\X)$ 
is the coproduct over the height $d$ points $s\in\X$ of the category 
$\cM\pc_\fin(\cA_s) = \bigcup_n \cM(A_s/p_s^n)$
of finite length $\cA_s$--sets.

Now $K(\cM\pc_\fin(\cA_s)) \simeq K(\G_s)$ by 
Theorem \ref{thm:devissage} (\ie by D\'evissage),
so we have long exact sequences
\begin{align}\label{eq:les} 
\cdots\to K_*\cM_{d+1}\pc(\X) \to &
K_*\cM_d\pc(\X) \to \bigoplus K_*(\G_s)\to \cdots
\\
\cdots\to K_*\cM_{d+1}(\X) \to &
K_*\cM_d(\X) \to \bigoplus G_*(\G_s)\to \cdots
\nonumber
\end{align}
The direct sums in these sequences 
are indexed by the points $s\in\X$ of codimension $d$.
Each of these sequences form exact couples, 
yielding fourth quadrant spectral sequences:
\begin{align}\label{eq:ss}
E_1^{p,q} = \bigoplus_{\height(s)=p} K_{-p-q}(\G_s) \Rightarrow K'_{-p-q}(\X)
\\
E_1^{p,q} = \bigoplus_{\height(s)=p} G_{-p-q}(\G_s)\Rightarrow G_{-p-q}(\X)
\nonumber
\end{align}
These are the analogues of the usual coniveau spectral sequence 
for the $K$-theory of coherent sheaves over a noetherian scheme;
see \cite[5.4]{Q} or \cite[V.9]{WK}.

\begin{ex}\label{DVM}
Let $D$ be a discrete valuation monoid with group of units $\G=D^\times$ and
parameter $t$, \ie $D=\G_+\wedge\N$. 
Since $D_0=(\G \times\Z)_+$, 
the localization sequence \cite[3.5]{HW} becomes
\[ \cdots  
K_n(\G) \map{0} K'_n(D) \to 
K_n(\G\times\Z) \map{d} K_{n-1}(U) \cdots
\]
(Recall that $K_n(\G)=K'_n(\G_+)$.)
The left map is zero by Additivity \cite[7.15]{CZ},
because for every $\G_+$-set $X$ there is a characteristic exact sequence
\[
   X\times\N \overset{t}\into X\times\N \onto X.
\]
Thus the $E_1$ page of the spectral sequence \eqref{eq:ss} for $K'(D)$
looks like:
\begin{equation*}
\vcenter{\xymatrix@C=1em@R=0em{
K_0(\G\times\Z) \ar[r] & 0       \\
K_1(\G\times\Z) \ar[r]^-d & K_0(\G) \ar[r] & 0\\
K_2(\G\times\Z) \ar[r] & K_1(\G) \ar[r] & 0 \ar[r] & 0\\
K_3(\G\times\Z) \ar[r] & K_2(\G) \ar[r] & 0 \ar[r] & 0
}}
\end{equation*}
Now $K_n(\G\times\Z) \cong K_n(\G) \oplus K_{n-1}(\G)$
by \ref{K(G)}, because $B\Z = S^1$,  and
\[ K(\G\times\Z) =
S^\infty B(\G\times\Z)_+ = S^\infty(S^1_+\wedge B\G_+)
 = S^1S^\infty(B\G)_+\vee S^\infty(B\G)_+.
\]
Hence $K_n(\G\times\Z)$
maps onto $K_{n-1}(\G)=\pi_{n-1}S^\infty(B\G)\vee \bS$
with kernel $K_n(\G)$.  Thus
$E_2$ is concentrated on the column $p=0$,
with $K'_q(D) = E_2^{0,-q}=K_q(\G)$.
\end{ex}

\begin{lem}\label{lem:d=1}
If $D$ is a discrete valuation monoid with units $\G$,\linebreak
$d_1^{0,-1}\colon K'_1(D_0)\to K_0(\G)=\Z$ is the surjection 
sending $[t]$ to $\pm1$.
\end{lem}

\begin{proof}
By Example \ref{DVM}, $K'_1(D_0)= (\G \times \Z) \oplus \pi_1^s$
and $K'_1(D) = \G \oplus \pi_1^s$ 
and $K'_1(D_0)/K'_1(D) \cong K_0(\G)=\Z$.
\end{proof}

\begin{subrem}
The analysis of the spectral sequence for $G_*(D)$ is more
complicated because the formula for $G_*(\G\times\Z)_+$ is more complicated.
\end{subrem}

\section{Weil divisors and $K'_0(\X)$}

\medskip
Recall that a monoid scheme is {\it normal} if each
$\cA_s$ is cancellative and if $a,b\in \cA_s$ satisfy
$(a/b)^n \in \cA_s$ then $a/b\in\cA_s$.
If $(\X,\cA_\X)$ is a normal (pc) monoid scheme 
then for each height~1 point $x$, $\cA_s$ is a 
Discrete Valuation monoid by \cite[2.6]{FW}.
The following definition is taken from \cite[4.2]{FW}.

\begin{defn}\label{def:Cl}
Let $(\X,\cA_\X)$ be a normal (pc) monoid scheme.
A {\it Weil divisor} is an element of $\Div(\X)$, 
the free abelian group on the height~1 points of $\X$.
If $a$ is a nonzero element of $\cA_0$, we can define the principal divisor
$\div(a)$ in $\Div(\X)$ by the same formula as in Algebraic Geometry.
The divisor class group $\Cl(\X)$ is the quotient of $\Div(\X)$
by the subgroup of principal divisors. Thus there is an exact sequence
\[
1 \to \cA_\X(X)^\times \to \cA_0^\times \map{\div} \Div(\X) \to \Cl(\X) \to 0.
\]
\end{defn}

\begin{subrem}
A monoid scheme $\X$ is {\it factorial} if and only if
$\Cl(\cA)=0$; if the stalks $\cA_s$ of $\cA$ are factorial then
$\Cl(\X)=\Pic(\X)$ by \cite[6.4]{FW}.
\end{subrem}

By Lemma \ref{lem:d=1}, the component
$K'_1(\cA_0) \map{d} K'_0(\cA_s)=\Z$ of the differential $d_1$ is the
projection $K'_1(\cA_0) \to \cA_0^\times$, followed by the map $\div$. Thus

\begin{thm}\label{thm:K'-Cl} 
If $\X$ is a normal (pc) monoid scheme, the differential
\linebreak 
$E_1^{0,-1} \map{d_1} E_1^{1,-1}$ is the principal divisor map $\div$,
and $E_\infty^{-1,-1}\cong \Cl(\X)$. 

Hence $K'_0(\X)$ surjects onto $\Z \oplus \Cl(\X)$.
\end{thm}

\begin{subrem}
If $A$ is a 2--dimensional normal monoid, the kernel of
$K'_0(A)\to\Z\oplus\Cl(X)$ is the subgroup generated by $A/\m$.
In general, the groups $E_2^{p,-p}(\X)$ are the 
higher class groups $W_p$ studied by Claborn and Fossum;
see \cite{CF} or \cite[Ex.\,II.6.4.3]{WK}.
\end{subrem}

The following example shows that the augmented rows of the spectral sequence 
aren't always exact when $A$ is normal.
\begin{ex}
Let $A$ be the (2-dimensional normal) submonoid of $\N^2$ generated 
by $x=(1,0)$, $z=(1,1)$ and $y=(1,2)$.
By \cite[4.2.1]{FW}, $\Cl(A)=\Z/2$ for $A$. Since $A/(y,z)\cong\N$,
and $A/(x,y,z)=*$, it follows that $E_2^{2,-2}=0$ and
$K'_0(A)\cong\Z\oplus \Z/2$.
\end{ex}

\section{Gersten's Conjecture for pc monoid schemes}
\label{sec:GC}

Quillen proved \cite[5.1]{Q} that if $R$ is a semilocal ring, smooth/$k$,
then the $K$-theory spectral sequence degenerates at $E_2$.
In this section, we give the analogue for monoids.

The following definition is taken from \cite[6.4]{CHWW};
the terminology reflects Theorem 6.5 in \cite{CHWW} 
that a (cancellative) monoid $A$ is 0-smooth if and only if
its monoid ring $k[A]$ is smooth over $k$ for every field $k$
of characteristic~0.

\begin{defn}\label{smooth}
A monoid $A$ is {\it 0-smooth} if it is a product $\G_+ \wedge (\N^n)$, 
where $\G$ is a finitely generated abelian group.
\end{defn}

\begin{thm}\label{thm:Gersten}
If $A$ is a 0-smooth monoid, the sequence
\[
0 \to K'_n(A) \to \bigoplus_{\height(s)=0} K'_n(A_s) \map{d}
\bigoplus_{\height(s)=1} K'_{n-1}(A_s) \map{d} \cdots 
\]
is exact for all $n$, where $d$ is the $E_1$ differential in 
the spectral sequence \eqref{eq:ss}.
The sequence remains exact if $K'_n$ is replaced by $G_n$.
\end{thm}

As in \cite[(5.9)]{Q}, this gives a flasque resolution of the sheaf
$\mathcal{K'}_n$ on any 0-smooth monoid scheme $\X$: 
\[
0\to \mathcal{K'}_n \to \bigoplus_{\height(s)=0}(i_s)_*K_n(\G_s) \to
\bigoplus_{\height(s)=1}(i_s)_*K_{n-1}(\G_s) \to \cdots
\]
and gives a canonical isomorphism 
\begin{equation}\label{eq:HpKq}
E_2^{pq}(\X) \cong H^p(\X,\mathcal{K}'_q).  
\end{equation}
Formula \eqref{eq:HpKq} remains valid
if $K'$-theory is replaced by $G$-theory.

\begin{proof}
We may suppose that $A=\G_+ \wedge \N^{m+1}$, 
with the maximal ideal generated by $x_0,\dots,x_m$. As noted by Quillen
\cite{Q}, it suffices to show that for every non-zerodivisor $t\in A$
that $\cM_d(A/tA)\to \cM_d(A)$ induces the zero map on $K$-theory.
Since $A$ is 0-smooth, $t=\gamma x_0^{n_0}\cdots x_m^{n_m}$ for
some $\gamma\in \G$, where some
$n_i>0$; we may suppose that $n_0>0$.

Set $B$ be the submonoid $\G_+\wedge\N^m$ on $x_1,...,x_m$.
Then $A$ is flat over $B$, so $A'=A\wedge_B(A/tA)$ is flat over $A/tA$,
and $A\to A/tA$ factors through $p:A\to A'$. There is a canonical
splitting $u:A'\to A/tA$, induced by $A\to A/tA$, so that
\[
\cM_d(A/tA) \map{u_*} \cM_d(A') \map{p^*} \cM_d(A)
\]
is the natural map. Now for any $A/tA$--set $X$, there is an
exact functor $v:\cM_d(A/tA) \map{} \cM_d(A')$ sending $X$ to
$v(X)=X\wedge\langle x_0\rangle$ 
and an exact sequence in $\cM_d(A')$:
\[
X\wedge\langle x_0\rangle \overset{x_0}\into 
X\wedge\langle x_0\rangle \onto u_*(X).
\]
By Additivity, $u_*\simeq 0$. Hence the map
$K(\cM_d(A/tA)) \map{} K(\cM_d(A))$
is zero, as asserted.
\end{proof}

\bigskip
\section{Regular localization}\label{sec:cat-thy}

The localization construction $\cM\to\cM/\cC$
generalizes from abelian categories
(due to Gabriel \cite{G}) to regular categories, including
$A\Sets$ and $A\Sets\pc$. In this section, we give a general framework 
into which our key examples fit. The following definition is taken from \cite[A.5]{BB}.

\begin{defn}\label{def:regular}
Let $\cM$ be a category with finite limits. An epimorphism
$Y\map{p} Z$ in $\cC$ is \emph{regular} if it is a coequalizer
$X\rightrightarrows Y\map{p} Z$.

We say that $\cM$ itself is \emph{regular} if regular epimorphisms are stable
under pullback by arbitrary morphisms and every \emph{kernel pair} admits a
coequalizer. Recall that the \emph{kernel pair} for a morphism $f\colon X\to Y$
is the diagram $X\times_Y X\rightrightarrows X$ given by the canonical maps in
the fiber product 
\begin{equation*}
\vcenter{\xymatrix@R=1em@C=1em{
X\times_Y X\ar[r]\ar[d]& X\ar[d]\\
X\ar[r]& Y
}}
\end{equation*}
The coequalizer of the kernel pair in a (concrete) regular category is the
quotient of $X$ given by identifying elements which become equal under $f$,
\ie the image of the morphism $f$. 

A sequence $X\smap{i} Y\smap{p} Z$ is called \emph{exact} if
$i=\ker (p)$ and $p=\coker (i)$. Kernels always exist in a regular
category, but cokernels do not necessarily exist; when they do, they are
regular epimorphisms and the coequalizer of their kernel pair. Thus given a
kernel $i\colon X\to Y$ we can always extend it to a short exact sequence.

A functor between regular categories is {\it regular} if it
preserves finite limits and regular epimorphisms (thus short exact sequences as well).
\end{defn}

The category of pointed sets is the prototype of a regular category; 
every surjection is a regular epimorphism (see \cite[A.5.6]{BB}).
In fact, the category of pointed objects in any elementary topos
is regular.

An abelian category is another (particularly nice) example of a regular
category. The regular epimorphisms are just cokernels. A regular functor
between abelian categories is just an exact functor.

\begin{rem}\label{reg=qexact}
Every regular category is quasi-exact, using the intrinsic notion of
short exact sequence in the regular category, and its associated
double category is a CGW category by \cite[2.3]{HW}.
All of the quasi-exact
categories considered in this paper are actually regular.
\end{rem}

One main feature of a regular category is that every morphism 
can be factored uniquely (up to isomorphism) as a regular epimorphism
(onto the image) followed by a monomorphism.

\begin{ex}
If $A$ is a (pointed) monoid, $A\Sets$ is regular, and
the forgetful functor $A\Sets\to\mathbf{Sets}$ is a regular functor
(because it has both a left and right adjoint; see \cite[2.5.1]{HW}).
If $A$ is a pc monoid, $A\Sets\pc$ is a regular subcategory of $A\Sets$.

Every monomorphism in $A\Sets$ (and hence $A\Sets\pc$) is a kernel, 
but not every epimorphism is a cokernel. In fact, 
the short exact sequences in $A\Sets$ and $A\Sets\pc$ have the form 
\[
X \overset{i}\into Y \overset{p}\onto Y/X.
\]
\end{ex}

The following definition is weaker than the related notions of a
{\it coherent} category and an {\it adhesive} category, 
but is just right for our purposes; see \cite[A.1]{J} for more details.

\begin{defn}\label{def:adherent}
A regular category $\cM$ is called \emph{adherent} 
if it satisfies the following axioms:
\begin{itemize}
\item[(A1)] Pushouts along kernels in $\cM$ 
are also pullbacks;
\item[(A2)] the image of every map is a kernel;
\item[(A3)] for each $X$ in $\cM$, the category of monics $X'\into X$
has both pushouts and pullbacks.
\end{itemize}
In particular, (A2) implies that every monomorphism is a kernel.
\end{defn}

Using De\,Morgan's laws on the underlying sets, it is easy to see that
$A\Sets$ and $A\Sets\pc$ are adherent categories. The 
category of pointed objects in any elementary topos, e.g., pointed
finite sets, is adherent; this gives a large class of examples.

\begin{construct}\label{pullyou}
Let $\cM$ be a adherent category, and
suppose that $X_i'\into X\onto X_i''$ ($i=1,2$) are
exact sequences in $\cM$.
Let $X'_{12}$ denote the pullback $X_1'\times_X X_2'$.
Write $X'$ for the pushout of $X'_1$ and $X'_2$ along $X'_{12}$,
which by abuse of notation we will also write $X'_1\cup X'_2$.
Note that $X'$ exists and $X'\into X$ is monic by axiom (A3).

Let $X''$ denote the cokernel of $X'\into X$
and write $X''_{12}$ for the pullback of
$X_1''$ and $X_2''$ along $X''$. Then we have exact sequences
\[
X'_{12} \into X \onto X''_{12} \quad\textrm{and}\quad
X' \into X \onto X''
\]
which fit into 
diagram \eqref{keydiagram}, where both squares are bicartesian.
\addtocounter{equation}{-1}
\begin{subequations}

\begin{equation}\label{keydiagram}
\vcenter{\xymatrix@R=1.9em{
X'_{12} \mono[r]\mono[d] & X'_1 \mono[d]\mono[dr] &&&\\
X'_2 \mono[r]  & X' \mono[r] & X \epi[r]\epi[dr]
                     & X''_{12} \epi[d]\epi[r] & X''_{2}\epi[d] \\
                          &&& X''_1 \epi[r] & X''
}}
\end{equation}
\end{subequations}
The top-left square is bicartesian by construction. The bottom-right square
is cartesian, so by \cite[A.1.4.3]{J} it is also cocartesian because cokernels are regular epimorphisms. Note
that $X'\into X\onto X''_{12}$ is \emph{not} claimed to be an exact sequence.

We obtain Diagram \eqref{keydiagram} as follows: 
by taking iterated pushouts along our string of monomorphisms, we obtain
\begin{equation*}
\vcenter{\xymatrix{
X_{12}'\mono[r]\ar[d]& X_1'\mono[r]\epi[d]& X'\mono[r]\epi[d]& X\epi[d]\\
\ast\mono[r]&X_1'/X_{12}'\mono[r]\ar[d]&X'/X_{12}'\mono[r]\epi[d]&X_{12}''\epi[d]\\
&\ast\ar[r]&X'/X_1'\mono[r]\ar[d]&X_1''\epi[d]\\
&&\ast\ar[r]&X'' .
}}
\end{equation*}
Note that all instances of $\twoheadrightarrow$ are cokernels. We can repeat this process with the chain
\begin{equation*}
\vcenter{\xymatrix{
X_{12}'\mono[r]& X_2'\mono[r]& X'\mono[r]& X
}}
\end{equation*}
to obtain the square
\begin{equation*}
\vcenter{\xymatrix{
X_{12}''\epi[r]\epi[d]&X_1''\epi[d]\\
X_2''\epi[r]&X'' .
}}
\end{equation*}
By construction, this is a pushout; as we have also identified it as a
pullback, this completes the construction. 
\end{construct}

\begin{defn}\label{defn:serresubcat}
  Let $\cM$ be an adherent category. A full subcategory
  $\cC\subset\cM$ is called \emph{Serre} if it is closed under finite
  limits and, for every short exact sequence $X\rightarrowtail
  Y\twoheadrightarrow Z$ in $\cM$, $Y$ is in $\cC$ if and only if
  $X,Z$ are both in $\cC$.
  
  Note that $\cC$ also forms an adherent category. It admits coequalizers of
  kernel pairs as, for any $f\colon X\to Y$ in $\cC$, the coequalizer
  of\linebreak $X\times_Y X\rightrightarrows Y$ is just $\im(f)$, which is a
  subobject of $Y$. Thus $\cC$ is regular, and axioms (A1), (A2), and (A3)
  hold because $\cC$ is a full subcategory closed under subobjects. 
\end{defn}

\begin{subrem}
The definition of a Serre subcategory of an abelian category
does not explicitly state that it be closed under finite limits. 
We need this assumption because we do not have biproducts; while there is a
canonical short exact sequence $X\rightarrowtail X\vee Z\twoheadrightarrow Z$
which makes $\cC$ closed under coproducts (and therefore pushouts, which are
quotients of coproducts), the product does not sit naturally in a short exact
sequence. To prove that the quotient category $\cM/\cC$ 
(to be defined shortly) has finite limits, 
we need $\cC$ to be closed under finite limits.
\end{subrem}

Before we define the quotient of an adherent category by a Serre subcategory,
we have the following technical lemma.

Let $X,Y$ be two objects in an adherent category $\cM$, and 
let $\cC$ be a Serre subcategory of $\cM$. 
Let $I=I_{X,Y}$ be the following category: 
its objects are the pairs $(X',Y'')$ 
for every pair of short exact sequences \linebreak
$X'\into X\onto X''$
and $Y'\into Y\onto  Y''$ such that $X'',Y'\in\cC$. 
The maps $(X_1',Y_1'')\to(X_2',Y_2'')$
are induced by composition with maps in $\cM$ of the form
$X_1'\into X_2'\into X$ and 
$Y\onto Y_1''\onto Y_2''$. Note that $I$ is equivalent to a partially-ordered set.

\begin{lem}\label{lem:filtcolimit}
Let $X,Y$ be two objects in an adherent category $\cM$, and 
let $\cC$ be a Serre subcategory of $\cM$. 
Then the category $I=I_{X,Y}$ is filtered, and 
$\Hom_{\cM}\colon I\to\mathbf{Sets}$ is a functor.
\end{lem}

\begin{proof}
To begin, fix two short exact sequences
$X'_i\into X \onto X''_i$ with $X''_i$ in $\cC$. 
Then the inclusions $X'_i\into X$ yield a span
\begin{equation*}
\vcenter{\xymatrix{
\Hom_\cM(X'_2,Y'') & \Hom_\cM(X,Y'')\ar[r]\ar[l]&\Hom_\cM(X'_1,Y'')
}}
\end{equation*}
\noindent Setting $X'_{12}=X_1'\times_X X_2''$,
we claim that $(X'_{12},Y'')$ is an upper bound of $(X'_1,Y'')$ and
$(X'_2,Y'')$ in $I$, i.e., that the cokernel $X''_{12}$ of
$X'_{12}\into X$ is in $\cC$.  By Construction~\ref{pullyou},
$X''_{12}$ is the pullback of $X_1''$ and $X_2''$ along $X''$. By
assumption, $X_1''$ and $X_2''$ are in $\cC$, and $X''$ is a quotient
of $X_2''$ (or $X_1''$) so is also in $\cC$. Thus $X''_{12}$ is a
(finite) limit of objects in $\cC$ and is also in this subcategory.

Dually, given short exact sequences
$Y'_i\into Y \onto Y''_i$ with $Y'_i$ in $\cC$,
we have the span
\begin{equation*}
\vcenter{\xymatrix{
\Hom_\cC(X,Y''_1) & \Hom_\cC(X,Y)\ar[r]\ar[l]&\Hom_\cC(X,Y''_2)
}}
\end{equation*}
Setting $Y''=Y_1''\vee_Y Y''_2$, 
it suffices to prove that the kernel $Y'$ of $Y\to Y''$ is in $\cC$.
Appealing to Construction~\ref{pullyou}, 
$Y'=Y_1'\cup Y_2'$, and
we have an identification of cokernels
\begin{equation*}
Y'/Y_1'=(Y_1'\cup Y_2')/Y_1'\cong Y_2'/(Y_1'\cap Y_2')=Y_2'/Y_{12}'.
\end{equation*} 
Since $Y_2'$ and hence $Y_2'/Y_{12}'$ is in $\cC$, this shows that
$Y'/Y_1'$ is in $\cC$.
Since $Y'_1$ is also in $\cC$, $Y'$ is in the Serre subcategory $\cC$,
as desired.
\end{proof}

\begin{thm}\label{thm:M/C-long}
Let $\cC$ be a Serre subcategory of an adherent category $\cM$.
Then we can define the quotient category $\cM/\cC$ in the following
way: the objects are the objects of $\cM$, and the morphisms are the sets
\begin{equation*}
\Hom_{\cM/\cC}(X,Y)=\colim_{I_{X,Y}} \Hom_\cM(X',Y'')
\end{equation*}
defined above. Then $\cM/\cC$ is a regular category, and $\cM\to\cM/\cC$ 
is a regular functor, initial amongst those sending $\cC$ to zero.
\end{thm}
\begin{proof}
We take inspiration from Gabriel's thesis \cite{G}, where this result is
proved for abelian categories. We first prove that our definition of
$\Hom_{\cM/\cC}$ behaves well with respect to composition. Suppose that
$f\in\Hom_{\cM/\cC}(X,Y)$ and $g\in\Hom_{\cM/\cC}(Y,Z)$. Then we can take
specific representatives $f\colon X'\to Y_1''$ and $g\colon Y_2'\to Z''$. To
compose these, we need to match the codomain of $f$ to the domain of $g$ by
taking equivalent representatives. 
Consider $Y_1'\rightarrowtail Y_1'\cup Y_2'$, 
where $Y_i'$ is the kernel of $Y\twoheadrightarrow Y_i''$ ($i=1,2$). 
Then the inverse image $X^*=f\inv(Y_1'\cup Y_2'/Y_1')$
is a subobject of $X'$ (and hence of $X$), fitting into the pullback diagram
\begin{equation*}
\vcenter{\xymatrix{
X^\ast\mono[r]\ar[d]_-{f^\ast}&X'\ar[d]^-f\\
(Y_1'\cup Y_2')/Y_1'\mono[r]&Y_1'' .
}}
\end{equation*}
The cokernel of $X^\ast\into X$ fits into a short exact sequence
\begin{equation*}
X'/X^\ast\rightarrowtail X/X^\ast\twoheadrightarrow X/X'=X''.
\end{equation*}
Cf.\,\cite[2.4]{HW}.
Because the source and quotient of this exact sequence are in $\cC$,
so is the middle; thus $X^\ast$ is a permissible domain in the colimit
defining $\Hom_{\cM/\cC}(X,Y)$. The codomain has kernel $Y'_1$, so it
is also permissible. The composition 
$f^\ast: X^\ast\to(Y_1'\cup Y_2')/Y_1'$ 
gives the same map as $f$ in the colimit.

For $g$, consider the pullback
$Y'_{12}=Y'_1\times_Y Y_2'$, and form the
quotient $Y'_{12} \into Y'_2 \onto Y'_2/Y'_{12}$;
taking the pushout gives a map
\begin{equation*}
\vcenter{\xymatrix{
Y_2'\ar[r]^-g\epi[d]&Z''\epi[d]\\
Y_2'/Y'_{12} \ar[r]^-{g^{\ast\ast}}&Z^{\ast\ast} .
}}
\end{equation*}
This pushout is constructed as follows: 
take the image factorization of the map
$Y'_{12} \into Z''$:
\begin{equation*}
\vcenter{\xymatrix{
Y'_{12} \mono[d]\ar@{-->}[r]&J\ar@{>-->}[d]\\
Y_2'\ar[r]^-g&Z''
}}
\end{equation*}
Because the image of every map in $\cM$ is a kernel by Axiom (A2), 
$J\rightarrowtail Z''$ is still a kernel. The pushout and the map 
$g^{\ast\ast}$ are created by taking cokernels vertically. 

By assumption, $Y_1'$ and $Y_2'$ are in $\cC$, 
so $Y'_{12}$ and $Y'_2/Y'_{12}$ are in $\cC$ as well.
Thus $g^{\ast\ast}$ is part of the colimit and represents the same map as
$g$. But by Construction \ref{pullyou} 
we have an identification 
\begin{equation*}
Y_1'\cup Y_2'/Y_1'\cong Y_2'/Y'_{12} 
\end{equation*}
which allows us to take the actual composition $g^{\ast\ast}\circ f^\ast$,
which gives the requisite composition law.
It is clear that composition is associative.

Having established how composition works, we have the following useful fact:
suppose that $X'\into X\onto X''$ is a short exact sequence with 
$X''$ in $\cC$. Then $X'\cong X$ in $\cM/\cC$.
Similarly, for any exact sequence $Y'\into Y\onto Y''$, if $Y'$ is in
$\cC$, then $Y\cong Y''$. Therefore given a map $f\colon X\to Y$ in $\cM/\cC$, 
the domain and codomain of the representative $X'\to Y''$ are 
isomorphic to the original ones.

To prove that the quotient category is regular, we first show that the
quotient is still finitely complete. Note that $\cM/\cC$ still has a zero
object, as the filtered colimits defining $\Hom_{\cM/\cC}(X,\ast)$ and
$\Hom_{\cM/\cC}(\ast,Y)$ are constant on the singleton set. Therefore the
existence of finite limits is equivalent to the existence of pullbacks.

Consider a cospan in $\cM/\cC$ of the form
\begin{equation*}
\vcenter{\xymatrix{
X_1\ar[r]^-{f_1} &Z& X_2\ar[l]_-{f_2}
}}
\end{equation*} 
To find the pullback in $\cM/\cC$,
we pick maps in $\cM$ that represent $f_1$ and $f_2$, 
and form a common codomain $Z''$. 
\goodbreak

\begin{equation*}
\vcenter{\xymatrix{
& Z_1\mono[d] &X'_2\ar[d]^-{f'_2} \\
Z_2 \mono[r]&Z\epi[r]\epi[d]&Z''_2\epi[d]\\
X_1' \ar[r]^-{f'_1}&Z_1''\epi[r]& Z''
}}
\end{equation*}
Here, as in  \eqref{keydiagram},
$Z''$ is the quotient of $Z$ by the pushout $Z'=Z'_1\cup Z'_2$ 
of the kernels $Z'_1$ and $Z'_2$,
along their intersection inside of $Z$; $Z'$ is still in $\cC$.
Thus $Z''$ is a common codomain.

As the total righthand and bottom compositions also represent the maps
$f_1,f_2$, we may pretend (up to renaming the $X'_i$ as $X_i$ and $Z''$ as $Z$) 
that this was our original situation. We know the fiber product 
$P=X_1\times_Z X_2$ exists in $\cM$
and we have projections $P\to X_1$ and $P\to X_2$ that
give rise to maps in $\cM/\cC$. We just need to prove that $P$ has the
requisite universal property.

Suppose that we have $T\in\cM/\cC$ fitting into a diagram
\begin{equation*}
\vcenter{\xymatrix{
T\ar@(r,u)[drr]\ar@(d,l)[ddr]\\
&P\ar[r]\ar[d]&X_2\ar[d]^-{f_2}\\
&X_1\ar[r]^-{f_1}&Z.
}}
\end{equation*} 
Then moving this picture over to $\cM$, we can pick representatives 
$h,k$ for the maps from $T$ to obtain 
\begin{equation*}
\vcenter{\xymatrix{
T'\ar@(r,u)[drrr]^-h\ar@(d,l)[dddr]_-k\\
&P\ar[r]\ar[d]&X_2\ar[d]^-{f_2}\epi[r]&X''_2\ar[dd]\\
&X_1\ar[r]^-{f_1}\epi[d]&Z\epi[dr]\\
&X_1''\ar[rr]&&Z'' .
}}
\end{equation*} 
There is a unique map from $T'$ to $P'':=X_1''\times_{Z''} X_2''$. 
Moreover, the kernel of $P\to P''$ is the pullback $X_1'\times_{Z'}X_2'$
of the kernels of $h$ and $k$.  Since $\cC$ is closed under finite limits,
$P''$ is an admissible quotient of $P$. Therefore the map 
$T'\to P''$ represents a map $T\to P$ in $\cM/\cC$. Note that any other
choice of representative for $h,k$ will lead to an equivalent map $T\to P$ in
the filtered colimit, which proves uniqueness. Therefore the quotient category
admits finite limits, and the quotient functor to $\cM/\cC$ preserves these, as
the limits are computed in $\cM$.

To see that  $\cM/\cC$ is regular, it remains to check that
the quotient still admits coequalizers of
kernel pairs which are stable under pullback. 
Suppose that $f\colon X\to Y$ is a map in $\cM/\cC$ and that we have already
replaced the domain and codomain so that $f\colon X\to Y$ is now in
$\cM$. Then we can take the coequalizer of its kernel pair in $\cM$ 
to obtain\linebreak $X\times_Y X\rightrightarrows Y\overset{p}\to Z$. 
We will show that $Y\map{p}Z$
represents the coequalizer in $\cM/\cC$.
By the previous argument, the kernel pair $X\times_Y X\rightrightarrows X$ 
in $\cM/\cC$ is computed in $\cM$, as it is a finite limit.

Suppose that $T\in\cM/\cC$ fits into the diagram
\begin{equation*}
\vcenter{\xymatrix{
X\times_Y X\ar[r]<0.5ex>\ar[r]<-0.5ex>&Y\ar[r]&T.
}}
\end{equation*}
Pick a representative $h\colon Y'\to T''$,
and take the pullback in $\cM$: 
\begin{equation*}
\vcenter{\xymatrix{
X'\mono[r]\ar[d]_{f'}&X\ar[d]^-f\\
Y'\mono[r]&Y.
}}
\end{equation*}
Now $X'\times_{Y'}X'\cong X\times_Y X$ in $\cM/\cC$, because
$Y'\into Y$ and hence $X'\into X$ are isomorphisms in $\cM/\cC$.
Thus the kernel pair for $f'\colon X'\to Y'$ in $\cM$ also
represents the kernel pair for $f$ in $\cM/\cC$.
In particular, the coequalizer
\begin{equation*}
X'\times_{Y'} X'\rightrightarrows Y'\overset{p'}\longrightarrow Z'
\end{equation*}
gives us a map $Z'\to T''$; because $Z'\cong Z$ we obtain a map $Z\to T$
in $\cM/\cC$ as required. Since this map is unique up to picking a different
representative for $Z\to T$, it is literally unique in
$\Hom_{\cM/\cC}(Z,T)$. Because we have already shown that finite limits and
coequalizers in $\cM/\cC$ are computed in $\cM$, 
coequalizers of kernel pairs are stable under pullback in $\cM/\cC$.

Finally, to show that the quotient functor $Q\colon\cM\to\cM/\cC$ is initial
amongst those sending $\cC$ to zero, suppose that $F\colon\cM\to\cN$ is
another regular functor that sends every object in $\cC$ to zero. Then for any
$X,Y\in\cM$, we obtain a map of filtered systems $I\to\mathbf{Sets}$ (where
$I=I_{X,Y}$ is the category of Lemma~\ref{lem:filtcolimit}) given by
\begin{equation*}
\Hom_\cM(X',Y'')\to\Hom_\cN(F(X'),F(Y''))
\end{equation*}
for each $(X',Y'')\in I$. Because both the cokernel of
$X'\rightarrowtail X$ and the kernel of
$Y\twoheadrightarrow Y''$ are in $\cC$,
we obtain isomorphisms $F(X')\cong F(X)$ and
$F(Y)\cong F(Y'')$, so that the codomain filtered system is essentially
constant with value $\Hom_\cN(F(X),F(Y))$. Taking the colimits yields a
natural map $\Hom_{\cM/\cC}(X,Y)\to \Hom_\cN(F(X),F(Y))$, which defines the
functor $\widetilde F\colon \cM/\cC\to\cN$ as required.
\end{proof}

We do not expect the quotient to be adherent in general, nor do we require it
for the application of \cite[8.6]{CZ}. 

\section{Comparison of localizations}\label{sec:compare}
We now prove the equivalence of the regular localization developed in the
previous section and the CGW-localization from \cite{CZ}.

\begin{lem}\label{lem:adherentACGW}
If $\cM$ is an adherent category, then, using the short exact
sequences of $\cM$, the associated double category $\cM\double$
is an ACGW-category.
\end{lem}

\begin{proof}
By Remark \ref{reg=qexact}, $\cM\double$ is a CGW category.
To add the letter `A', we need to check a few more axioms, listed in
\cite[5.5--5.6]{CZ}. Axiom (P) is the assertion that
monomorphisms are closed under pullbacks, and
cokernels are closed under pushouts (as we showed in the proof of
Lemma~\ref{lem:filtcolimit}). 

For axiom (U), we need to check that commutative squares of
mono-morphisms 
give rise (upon taking cokernels) to pullback squares of cokernels;
this was again proven in Construction~\ref{pullyou} in the specific
case of $X'=X_1'\cup X_2'$, but the proof in the general case is
identical. The mixed pullback square is defined using the
factorizations of morphisms in a regular category; the object
$X\oslash_Y Z$, defined in \cite[5.6]{CZ}
for a composable pair $X\rightarrowtail Y\twoheadrightarrow Z$,
is just the image of the composition. The compatibility condition
follows from the uniqueness of such factorizations.

Axiom (S) concerns the pullback of monomorphisms and their associated
cokernels. We have already checked everything required for this axiom
in Construction~\ref{pullyou}. In particular, the restricted pushout of this
axiom is just the pushout in the ambient category $\cM$.

Finally, axiom (PP) uses the full strength of an adherent
category. Because restricted pushouts for us are pushouts, Axiom (A3)
guarantees that these always exist always exist. Pushouts also exist
along cokernels, as we demonstrated in Construction~\ref{pullyou}, and
there is no further compatibility to check as restricted pushouts are
just pushouts (so they are appropriately functorial).
\end{proof}

\begin{lem}\label{lem:M/C[d]}
Let $\cC$ be a Serre subcategory of an adherent category $\cM$.
Then the CGW-category $(\cM/\cC)\double$ is equivalent to the 
CGW-quotient $\cM\double\backslash\cC\double$.
\end{lem}

As noted in Section~\ref{sec:cgw}, it follows that 
$K(\cM/\cC)\cong K(\cM\double\backslash\cC\double)$.

\begin{proof}
The definition of the CGW-quotient can be found at \cite[8.1]{CZ}.

The m-morphisms $V\rightarrowtail Z$ of the CGW-quotient 
$\cM\double\backslash \cC\double$ are defined as compositions
\begin{equation*}
\vcenter{\xymatrix{
V\epi[r]|\bullet & W & X\mono[l]|\bullet & Y\epi[l]|\bullet\mono[r]& Z
}}
\end{equation*}
where the decoration $\bullet$ denotes that the (co)kernel of the map is in
the Serre subcategory $\cC$. The monomorphisms $V\rightarrowtail Z$ of
$\cM/\cC$ are those that have a representative $V'\to Z''$ which has its
kernel in $\cC$. Put another way, they are defined by
\begin{equation*}
\vcenter{\xymatrix{
V& V'\mono[l]|\bullet\mono[r] & Z'' & Z\epi[l]|\bullet.
}}
\end{equation*}
We can take the pushout $Y$ of the span of monomorphisms to obtain an
equivalent expression 
\begin{equation*}
\vcenter{\xymatrix{
V \mono[r]& Y & Z'' \mono[l]|\bullet& Z\epi[l]|\bullet.
}}
\end{equation*}
By inspection, 
kernel and the cokernel of $Z\to Y$ belong to $\cC$. This means that any
monomorphism $V\to Z$ can be described as a right fraction in a way identical
to \cite[8.3]{CZ}. The same reasoning goes through to identify the
m-morphisms in $\cM/\cC$ with those in $\cM\double\backslash\cC\double$; the
argument for the e-morphisms follows by dual reasoning.
\end{proof}

We will make use of the following lemma concerning isomorphisms in the
quotient category $\cM/\cC$. 

\begin{lem}\label{lem:monorep}
Let $\cC$ be a Serre subcategory of an adherent category $\cM$. Let $h\colon
X\dashrightarrow V$ be an isomorphism in $\cM/\cC$. Then there exists a
representative $h\colon X'\to V''$ of $h$ in $\cM$ which is a retract, and
hence a monomorphism.
\end{lem}

\begin{proof}
Because $h$ is an isomorphism, we know there must exist some inverse 
$p\colon V\dashrightarrow X$ such that $p\circ h=\id_X$ in $\cM/\cC$. 
If we pick representatives in $\cM$ for $h,p$ and compose them as in the proof
of Theorem~\ref{thm:M/C-long}, we obtain
\begin{equation*}
\vcenter{\xymatrix{
X'\ar[r]^{h}\ar@(dr,dl)[rr]_-=&V''\ar[r]^-p&X'
}}
\end{equation*}
for some $X'\into X$ and $V\onto V''$. 
Because $h$ is a retract, it is a monomorphism.
\end{proof}

\begin{thm}\label{thm:MC=MC-long}
Let $\cC$ be a Serre subcategory of an adherent category $\cM$.  
Then the double categories $\cC\double$ and $\cM\double$ satisfy
conditions (W) and (CGW) of \cite[8.6]{CZ}. 
We therefore obtain a homotopy fiber sequence on K-theory:
\begin{equation*}
K(\cC)\to K(\cM)\to K(\cM/\cC).
\end{equation*}
\end{thm}

Special cases of Theorem \ref{thm:MC=MC-long} are given in
 \cite[8.3]{CZ} and \cite[3.5]{HW}.

\begin{subrem}\label{rem:(E)}
Condition (E) of \cite[8.6]{CZ} is superfluous. It 
is used only in the proof of \cite[10.22]{CZ}, a lemma in service of their
main localization theorem. However, in the absence of condition (E),
the {\it isomorphism} of categories in \emph{loc.~cit.}~may be weakened to
an {\it equivalence} of categories, which still induces the requisite
homotopy equivalence for the main theorem. We do not know whether
condition (E) holds in general for the CGW-categories associated to a
Serre subcategory of an adherent category, but fortunately it is
extraneous.
\end{subrem}

\begin{proof}
  We can apply \cite[8.6]{CZ} to obtain the homotopy fiber sequence as
  soon as we verify these conditions. Condition (CGW) holds by
  Theorem~\ref{thm:M/C-long} because $\cM/\cC$ is a regular category
  which is the ambient category of $\cM\double\backslash\cC\double$.

  For condition (W), we will prove that the categories $\mathcal
  I_V^m$ and $\mathcal I_V^e$ are filtered for all $V\in\cM$. The
  category $\mathcal I_V^m$ has objects $(X,\phi)$ with $\phi\colon
  X\dashrightarrow V$ an isomorphism in $\cM/\cC$; the maps of
  $\mathcal I_V^m$ (under our hypotheses) are maps $g\colon X_1\to
  X_2$ in $\cM$ which are sent to isomorphisms over $V$ in
  $\cM/\cC$. The category $\mathcal I_V^m$ is isomorphic to $I_V^e$ by
  appealing to image factorizations, just as in the abelian case,
  cf.~\cite[8.7]{CZ}. Therefore it suffices to prove that $I_V^m$ is
  filtered.

  Recall that a nonempty category $I$ is filtered if two conditions
  hold. First, for any two objects $X,Y\in I$, there exists a third
  object $Z$ and maps $X\to Z$ and $Y\to Z$, \ie $Z$ is an upper bound
  for $X$ and $Y$. Second, any two parallel arrows $g_1,g_2\colon X\to
  Y$ admit a weak coequalizer, \ie some $h\colon Y\to Z$ such that
  $h\circ g_1=h\circ g_2$. These conditions may be combined by saying
  that any finite diagram in $I$ admits a (non-unique) cocone.

  Suppose we have two objects $(X_1,\phi_1)$ and $(X_2,\phi_2)$ of
  $I^m_V$. If we choose representatives of $\phi_1$ and $\phi_2$ using
  Lemma~\ref{lem:monorep}, up to identifying their codomains we obtain
\begin{equation*}
\vcenter{\xymatrix{
&X_1'\mono[dl]|-\bullet\mono[r]^-{\phi_1}&V''&X'_2\mono[l]_-{\phi_2}\mono[dr]|-\bullet\\
X_1&&V\epi[u]|-\bullet&&X_2
}}
\end{equation*}
where the vertical arrows are maps in $\cM$ which become isomorphisms
in $\cM/\cC$ and the horizontal arrows are monomorphisms in
$\cM$. Taking pushouts on the left and right we obtain
\begin{equation*}
\vcenter{\xymatrix{
&X_1'\mono[dl]|-\bullet\mono[r]^-{\phi_1}&V''\mono[dl]|-\bullet\mono[dr]|-\bullet&X'_2\mono[l]_-{\phi_2}\mono[dr]|-\bullet\\
X_1\ar[r]&Y_1&V\epi[u]|-\bullet&Y_2&X_2\ar[l]
}}
\end{equation*}
and finally pushing out the central span gives
\begin{equation*}
\vcenter{\xymatrix{
&X_1'\mono[dl]|-\bullet\mono[r]^-{\phi_1}&V''\mono[dl]|-\bullet\mono[dr]|-\bullet&X'_2\mono[l]_-{\phi_2}\mono[dr]|-\bullet\\
X_1\ar[r]&Y_1\mono[dr]|-\bullet&V\epi[u]|-\bullet&Y_2\mono[dl]|-\bullet&X_2\ar[l]\\
&&Y\ar@{-->}[u]^-\psi
}}
\end{equation*}
The map $\psi$ does not exist until we pass to the quotient category. The object $(Y,\psi)$ is then an upper bound for $(X_1,\phi_1)$ and $(X_2,\phi_2)$, which
proves the first condition of $I_V^m$ being filtered.

Supposing now that we have two parallel morphisms in $I_V^m$
\begin{equation*}
g_1,g_2\colon (X,\phi)\to(Y,\psi),
\end{equation*}
we need to find a weak coequalizer. Since these are maps over $V$, we have an
equality of isomorphisms $\psi\circ g_1=\psi\circ g_2=\phi$ in $\cM/\cC$. 
Therefore $\psi$ is nearly the map we want, but it does not exist in $\cM$. 

To remedy this, we consider the map 
$h=\psi^{-1}\circ\phi\colon X\dashrightarrow Y$ in $\cM/\cC$. 
By Lemma~\ref{lem:monorep}, we may take a monic representative 
$h\colon X'\to Y''$ with section $p\colon  Y''\to X'$ 
such that $p\circ h=\id_{X'}$ in $\cM$. Since $p$ represents
$\phi^{-1}\circ\psi$, we obtain a weak coequalizer by composing $p$ with
$g_i$: 
\begin{equation*}
\vcenter{\xymatrix{
X\ar[r]<-.5ex>_-{g_2}\ar[r]<.5ex>^-{g_1}&Y\epi[r]&Y''\ar[r]^-p&X'.
}}
\end{equation*}
Specifically, the weak coequalizer is $(Y,\psi)\to (X',\phi|_{X'})$ given by
the canonical projection to $Y''$ followed by $p$. 

We conclude that $I_V^m$ and $I_V^e$ are filtered, thus $\cC$ is both m- and
e-well represented in $\cM$. This proves that $\cC\subset\cM$ satisfies
condition (W), so applying \cite[8.6]{CZ} we complete the proof. 
\end{proof}

\begin{rem}
  A recent preprint by Sarazola and Shapiro in \cite{SS} also studies
  $\cM\backslash\cC$ using an alternative approach to ours which
  focuses more on the double-categorical aspects of
  Campbell-Zakharevich's original work \cite{CZ}.
\end{rem}


\newpage
\begin{thebibliography}{10}

\bibitem[Bar]{Bar}
C. Barwick,
Multiplicative structures on algebraic {K}-theory,
Doc. Math. 20 (2015), 859–878.

\bibitem[BB]{BB}
F, Borceux and D. Bourn,
Mal'cev, protomodular, homological and semiabelian categories,
Math. Appl. 566, Kluwer, 2004.

\bibitem[CZ]{CZ}
J. Campbell and I. Zakharevich,
Devissage and localization for the Grothendieck spectrum of varieties,
preprint (2018), arXiv:1811.08014v2.

\bibitem[CDD]{CDD}
G. Carlsson, C. Douglas and B. Dundas,
Higher topological cyclic homology and the Segal conjecture for tori,
Adv. Math. 226 (2011), 1823--1874. 

\bibitem[CLS]{CLS}
C. Chu, O. Lorscheid and R. Santhanam,
Sheaves and $K$-theory for $\mathbb{F}_1$-schemes,
Adv. Math. 229 (2012), 2239--2286.

\bibitem[CF]{CF}
L. Claborn and R. Fossum,
Generalizations of the notion of class group,
Illinois J. Math. 12 (1968), 228--253.

\bibitem[CHWW]{CHWW}
G. Corti\~nas, C. Haesemayer, M.\,E. Walker and C. Weibel,
Toric varieties, monoid schemes and cdh descent,
J. Reine Angew. Math. 698 (2015), 1--54. 

\bibitem[D]{Deitmar}
A. Deitmar, 
Remarks on zeta functions and $K$-theory over $\F_1$,
Proc. Japan Acad. Ser. A Math. Sci. 82 (2006), 141--146. 

\bibitem[FW]{FW}
J. Flores and C. Weibel, 
Picard groups and class groups of monoid schemes, 
J. Algebra 415 (2014), 247--263.

\bibitem[G]{G} 
P. Gabriel,
Des cat\'egories ab\'eliennes,
Bull. SMF 90 (1962), 323--448.

\bibitem[HW]{HW}
C. Haesemeyer and C. Weibel,
The $K'$--theory of monoid sets,
Proc. AMS, to appear.

\bibitem[J]{J}
P. Johnstone,
Sketches of an elephant: a topos theory compendium,
vol. 43 of Oxford Logic Guides, 
Oxford University Press, New York, 2002

\bibitem[Q]{Q}
D. Quillen,
Higher algebraic $K$-theory I,
pp. 85--147 in Lecture Notes in Math. 341, Springer, 1973.

\bibitem[SS]{SS}
M. Sarazola and B. Shapiro,
A Gillet-Waldhausen Theorem for chain complexes of sets,
arXiv:2107.07701

\bibitem[T]{T}
P. Taylor,
Practical foundations of mathematics,
vol. 59 of Cambridge Studies in Advanced Mathematics,
Cambridge University Press, Cambridge, 1999.


\bibitem[WK]{WK}
C. Weibel,
{\it The $K$-book}, 
Grad. Studies in Math. 145, AMS, 2013.

\end{thebibliography}
\end{document}